\documentclass[11pt]{article}
\usepackage{mathrsfs}
\usepackage{amscd}
\usepackage{amsmath,amsfonts,amssymb,amscd}
\usepackage{indentfirst,graphics,epsfig,psfrag}
\input{epsf}
\usepackage{ifpdf}
\usepackage{enumerate}
\usepackage{appendix}
\usepackage{enumerate}

\usepackage{lineno}

\makeatletter \@addtoreset{figure}{section} \makeatother
\makeatletter
\long\def\@makecaption#1#2{%
   \vskip 10\p@
   \setbox\@tempboxa\hbox{{#1}\ \ #2}%
   \ifdim \wd\@tempboxa >\hsize

       {#1}\ \ #2\par
   \else
       \hbox to\hsize{\hfil\box\@tempboxa\hfil}%
   \fi}
\makeatother

\newtheorem{thm}{Theorem}

\newtheorem{lem}{Lemma}

\newtheorem{obs}{Observation}
\newtheorem{pro}{Proposition}

\newcommand{\qed}{{\hfill\rule{3pt}{7pt}}}

\def\qed{\hfill \rule{4pt}{7pt}}

\begin{document}
\title{\textbf{The vertex-rainbow index of a graph}
\footnote{Supported by the National Science Foundation of China (No.
11161037) and the Science Found of Qinghai Province (No.
2014-ZJ-907).}}
\author{
\small  Yaping Mao\footnote{E-mail: maoyaping@ymail.com}\\[0.2cm]
\small Department of Mathematics, Qinghai Normal\\
\small University, Xining, Qinghai 810008, China\\
\small }
\date{}
\maketitle
\begin{abstract}
The $k$-rainbow index $rx_k(G)$ of a connected graph $G$ was
introduced by Chartrand, Okamoto and Zhang in 2010. As a natural
counterpart of the $k$-rainbow index, we introduced the concept of
$k$-vertex-rainbow index $rvx_k(G)$ in this paper. For a graph
$G=(V,E)$ and a set $S\subseteq V$ of at least two vertices,
\emph{an $S$-Steiner tree} or \emph{a Steiner tree connecting $S$}
(or simply, \emph{an $S$-tree}) is a such subgraph $T=(V',E')$ of
$G$ that is a tree with $S\subseteq V'$. For $S\subseteq V(G)$ and
$|S|\geq 2$, an $S$-Steiner tree $T$ is said to be a
\emph{vertex-rainbow $S$-tree} if the vertices of $V(T)\setminus S$
have distinct colors. For a fixed integer $k$ with $2\leq k\leq n$,
the vertex-coloring $c$ of $G$ is called a \emph{$k$-vertex-rainbow
coloring} if for every $k$-subset $S$ of $V(G)$ there exists a
vertex-rainbow $S$-tree. In this case, $G$ is called
\emph{vertex-rainbow $k$-tree-connected}. The minimum number of
colors that are needed in a $k$-vertex-rainbow coloring of $G$ is
called the \emph{$k$-vertex-rainbow index} of $G$, denoted by
$rvx_k(G)$. When $k=2$, $rvx_2(G)$ is nothing new but the
vertex-rainbow connection number $rvc(G)$ of $G$. In this paper,
sharp upper and lower bounds of $srvx_k(G)$ are given for a
connected graph $G$ of order $n$,\ that is, $0\leq srvx_k(G)\leq
n-2$. We obtain the Nordhaus-Guddum results for $3$-vertex-rainbow
index, and show that $rvx_3(G)+rvx_3(\overline{G})=4$ for $n=4$ and
$2\leq rvx_3(G)+rvx_3(\overline{G})\leq n-1$ for $n\geq 5$. Let
$t(n,k,\ell)$ denote the minimal size of a connected graph $G$ of
order $n$ with $rvx_k(G)\leq \ell$, where $2\leq \ell\leq n-2$ and
$2\leq k\leq n$.
The upper and lower bounds for $t(n,k,\ell)$ are also obtained. \\[2mm]
\end{abstract}

\section{Introduction}

The rainbow connections of a graph which are applied to measure the
safety of a network are introduced by Chartrand, Johns, McKeon and
Zhang \cite{Chartrand}. Readers can see \cite{Chartrand, Chartrand4,
Chartrand3} for details. Consider an edge-coloring (not necessarily
proper) of a graph $G=(V,E)$. We say that a path of $G$ is
\emph{rainbow}, if no two edges on the path have the same color. An
edge-colored graph $G$ is \emph{rainbow connected} if every two
vertices are connected by a rainbow path. The minimum number of
colors required to rainbow color a graph $G$ is called \emph{the
rainbow connection number}, denoted by $rc(G)$. In
\cite{M.Krivelevich}, Krivelevich and Yuster proposed a similar
concept, the concept of vertex-rainbow connection. A vertex-colored
graph $G$ is \emph{vertex-rainbow connected} if every two vertices
are connected by a path whose internal vertices have distinct
colors, and such a path is called a \emph{vertex-rainbow path}. The
\emph{vertex-rainbow connection number} of a connected graph $G$,
denoted by $rvc(G)$, is the smallest number of colors that are
needed in order to make $G$ vertex-rainbow connected. For more
results on the rainbow connection and vertex-rainbow connection, we
refer to the survey paper \cite{LiSun} of Li, Shi and Sun and a new
book \cite{LiSun1} of Li and Sun. All graphs considered in this
paper are finite, undirected and simple. We follow the notation and
terminology of Bondy and Murty \cite{Bondy}, unless otherwise
stated.

For a graph $G=(V,E)$ and a set $S\subseteq V$ of at least two
vertices, \emph{an $S$-Steiner tree} or \emph{a Steiner tree
connecting $S$} (or simply, \emph{an $S$-tree}) is a such subgraph
$T=(V',E')$ of $G$ that is a tree with $S\subseteq V'$. A tree $T$
in $G$ is a \emph{rainbow tree} if no two edges of $T$ are colored
the same. For $S\subseteq V(G)$, a \emph{rainbow $S$-Steiner tree}
(or simply, \emph{rainbow $S$-tree}) is a rainbow tree connecting
$S$. For a fixed integer $k$ with $2\leq k\leq n$, the edge-coloring
$c$ of $G$ is called a \emph{$k$-rainbow coloring} if for every
$k$-subset $S$ of $V(G)$ there exists a rainbow $S$-tree. In this
case, $G$ is called \emph{rainbow $k$-tree-connected}. The minimum
number of colors that are needed in a $k$-rainbow coloring of $G$ is
called the \emph{$k$-rainbow index} of $G$, denoted by $rx_k(G)$.
When $k=2$, $rx_2(G)$ is the rainbow connection number $rc(G)$ of
$G$. For more details on $k$-rainbow index, we refer to \cite{CLS,
CLS2, Chartrand2, CLYZ, LSYZ, LSYZ2}.

Chartrand, Okamoto and Zhang \cite{Chartrand3} obtained the
following result.

\begin{thm}{\upshape\cite{Chartrand2}}\label{th1-1}
For every integer $n\geq 6$, $rx_3(K_n)=3$.
\end{thm}

As a natural counterpart of the $k$-rainbow index, we introduce the
concept of $k$-vertex-rainbow index $rvx_k(G)$ in this paper. For
$S\subseteq V(G)$ and $|S|\geq 2$, an $S$-Steiner tree $T$ is said
to be a \emph{vertex-rainbow $S$-tree} or \emph{vertex-rainbow tree
connecting $S$} if the vertices of $V(T)\setminus S$ have distinct
colors. For a fixed integer $k$ with $2\leq k\leq n$, the
vertex-coloring $c$ of $G$ is called a \emph{$k$-vertex-rainbow
coloring} if for every $k$-subset $S$ of $V(G)$ there exists a
vertex-rainbow $S$-tree. In this case, $G$ is called
\emph{vertex-rainbow $k$-tree-connected}. The minimum number of
colors that are needed in a $k$-vertex-rainbow coloring of $G$ is
called the \emph{$k$-vertex-rainbow index} of $G$, denoted by
$rvx_k(G)$. When $k=2$, $rvx_2(G)$ is nothing new but the
vertex-rainbow connection number $rvc(G)$ of $G$. It follows, for
every nontrivial connected graph $G$ of order $n$, that
$$
rvx_2(G)\leq rvx_3(G)\leq \cdots \leq rvx_n(G).
$$

Let $G$ be the graph of Figure 1 $(a)$. We give a vertex-coloring
$c$ of the graph $G$ shown in Figure 1 $(b)$. If $S=\{v_1,v_2,v_3\}$
(see Figure 1 $(c)$), then the tree $T$ induced by the edges in
$\{v_1u_1,v_2u_1,u_1u_4,u_4v_3\}$ is a vertex-rainbow $S$-tree. If
$S=\{u_1,u_2,v_3\}$, then the tree $T$ induced by the edges in
$\{u_1u_2,u_2u_4,u_4v_3\}$ is a vertex-rainbow $S$-tree. One can
easily check that there is a vertex-rainbow $S$-tree for any
$S\subseteq V(G)$ and $|S|=3$. Therefore, the vertex-coloring $c$ of
$G$ is a $3$-vertex-rainbow coloring. Thus $G$ is vertex-rainbow
$3$-tree-connected.
\begin{figure}[!hbpt]
\begin{center}
\includegraphics[scale=0.85]{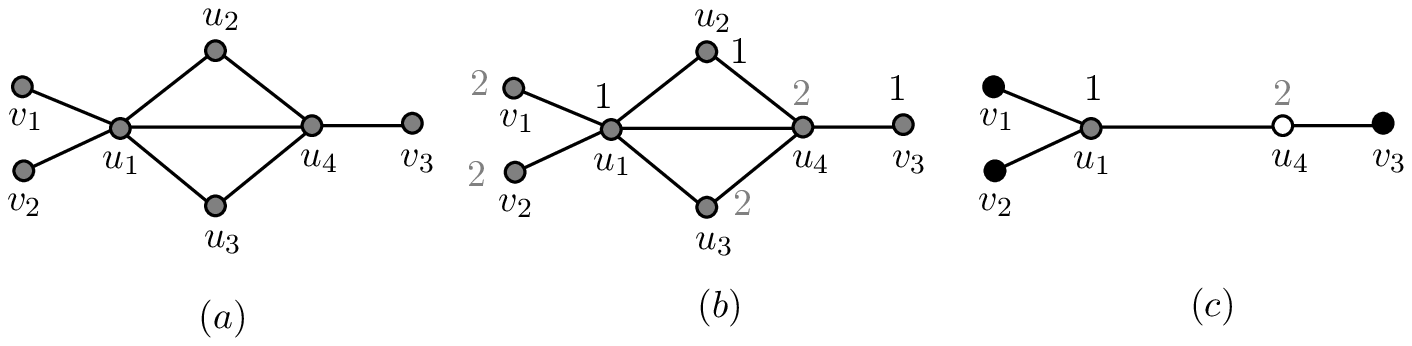}\\
Figure 1: Graphs for the basic definitions.
\end{center}
\end{figure}

In some cases $rvx_k(G)$ may be much smaller than $rx_k(G)$. For
example, $rvx_k(K_{1,n-1})=1$ while $rx_k(K_{1,n-1})=n-1$ where
$2\leq k\leq n$. On the other hand, in some other cases, $rx_k(G)$
may be much smaller than $rvx_k(G)$. For $k=3$, we take $n$
vertex-disjoint cliques of order $4$ and, by designating a vertex
from each of them, add a complete graph on the designated vertices.
This graph $G$ has $n$ cut-vertices and hence $rvx_3(G)\geq n$. In
fact, $rvx_3(G)=n$ by coloring only the cut-vertices with distinct
colors. On the other hand, from Theorem \ref{th1-1}, it is not
difficult to see that $rx_3(G)\leq 9$. Just color the edges of the
$K_n$ with, say, color $1,2,3$ and color the edges of each clique
with the colors $4,5,\cdots,9$.

Steiner tree is used in computer communication networks (see
\cite{Du}) and optical wireless communication networks (see
\cite{Cheng}). As a natural combinatorial concept, the rainbow index
and the vertex-rainbow index can also find applications in
networking. Suppose we want to route messages in a cellular network
in such a way that each link on the route between more than two
vertices is assigned with a distinct channel. The minimum number of
channels that we have to use is exactly the rainbow index and
vertex-rainbow index of the underlying graph.

The Steiner distance of a graph, introduced by Chartrand,
Oellermann, Tian and Zou \cite{Chartrand2} in 1989, is a natural
generalization of the concept of classical graph distance. Let $G$
be a connected graph of order at least $2$ and let $S$ be a nonempty
set of vertices of $G$. Then the \emph{Steiner distance} $d(S)$
among the vertices of $S$ (or simply the distance of $S$) is the
minimum size among all connected subgraphs whose vertex sets contain
$S$. Let $n$ and $k$ be two integers with $2\leq k\leq n$. The
\emph{Steiner $k$-eccentricity $e_k(v)$} of a vertex $v$ of $G$ is
defined by $e_k(v)=\max \{d(S)\,|\,S\subseteq V(G), |S|=k,~and~v\in
S \}$. The \emph{Steiner $k$-diameter} of $G$ is $sdiam_k(G)=\max
\{e_k(v)\,|\,v\in V(G)\}$. Clearly, $sdiam_k(G)\geq k-1$.

Then, it is easy to see the following results.

\begin{pro}\label{pro1}
Let $G$ be a nontrivial connected graph of order $n$. Then
$rvx_k(G)=0$ if and only if $sdiam_k(G)=k-1$.
\end{pro}

\begin{pro}\label{pro2}
Let $G$ be a nontrivial connected graph of order $n \ (n\geq 5)$,
and let $k$ be an integer with $2\leq k\leq n$. Then
$$
0\leq rvx_k(G)\leq n-2.
$$
\end{pro}
\begin{pf}
We only need to show $rvx_k(G)\leq n-2$. Since $G$ is connected,
there exists a spanning tree of $G$, say $T$. We give the internal
vertices of the tree $T$ different colors. Since $T$ has at most two
leaves, we must use at most $n-2$ colors to color all the internal
vertices of the tree $T$. Color the leaves of the tree $T$ with the
used colors arbitrarily. Note that such a vertex-coloring makes $T$
vertex-rainbow $k$-tree-connected. Then $rvx_k(T)\leq n-2$ and hence
$rvx_k(G)\leq rvx_k(T)\leq n-2$, as desired. \qed
\end{pf}

\begin{obs}\label{obs1}
Let $K_{s,t}$, $K_{n_1,n_2,\ldots,n_k}$, $W_{n}$ and $P_n$ denote
the complete bipartite graph, complete multipartite graph, wheel and
path, respectively. Then

$(1)$ For integers $s$ and $t$ with $s\geq 2,t \geq 1$,
$rvc(K_{s,t})=1$.

$(2)$ For $k\geq 3$, $rvx_k(K_{n_1,n_2,\ldots,n_k})=1$.

$(3)$ For $n\geq 4$, $rvx_k(W_{n})=1$.

$(4)$ For $n\geq 3$, $rvx_k(P_n)=n-2$.
\end{obs}

Let $\mathcal {G}(n)$ denote the class of simple graphs of order $n$
and $\mathcal {G}(n,m)$ the subclass of $\mathcal {G}(n)$ having
graphs with $n$ vertices and $m$ edges. Give a graph parameter
$f(G)$ and a positive integer $n$, the \emph{Nordhaus-Gaddum
(\textbf{N-G}) Problem} is to determine sharp bounds for: $(1)$
$f(G)+f(\overline{G})$ and $(2)$ $f(G)\cdot f(\overline{G})$, as $G$
ranges over the class $\mathcal {G}(n)$, and characterize the
extremal graphs. The Nordhaus-Gaddum type relations have received
wide attention; see a recent survey paper \cite{Aouchiche} by
Aouchiche and Hansen.

Chen, Li and Lian \cite{CLLian} gave sharp lower and upper bounds of
$rx_k(G)+rx_k(\overline{G})$ for $k=2$. In \cite{CLLiu}, Chen, Li
and Liu obtained sharp lower and upper bounds of
$rvx_k(G)+rvx_k(\overline{G})$ for $k=2$. In Section $2$, we
investigate the case $k=3$ and give lower and upper bounds of
$rvx_3(G)+rvx_3(\overline{G})$.

\begin{thm}\label{th2}
Let $G$ and $\overline{G}$ be a nontrivial connected graph of order
$n$. If $n=4$, then $rvx_3(G)+rvx_3(\overline{G})=4$. If $n\geq 5$,
then we have
$$
2\leq rvx_3(G)+rvx_3(\overline{G})\leq n-1.
$$
Moreover, the bounds are sharp.
\end{thm}

Let $s(n,k,\ell)$ denote the minimal size of a connected graph $G$
of order $n$ with $rx_k(G)\leq \ell$, where $2\leq \ell\leq n-1$ and
$2\leq k\leq n$. Schiermeyer \cite{Schiermeyer} focused on the case
$k=2$ and gave exact values and upper bounds for $s(n,2,\ell)$.
Later, Li, Li, Sun and Zhao \cite{LLSZ} improved Schiermeyer's lower
bound of $s(n,2,2)$ and get a lower bound of $s(n,2,\ell)$ for
$3\leq \ell \leq \lceil\frac{n}{2}\rceil$.

In Section $3$, we study the vertex case. Let $t(n,k,\ell)$ denote
the minimal size of a connected graph $G$ of order $n$ with
$rvx_k(G)\leq \ell$, where $2\leq \ell\leq n-2$ and $2\leq k\leq n$.
We obtain the following result in Section $3$.

\begin{thm}\label{th3}
Let $k,n,\ell$ be three integers with $2\leq \ell\leq n-3$ and
$2\leq k\leq n$. If $k$ and $\ell$ has the different parity, then
$$
n-1\leq t(n,k,\ell)\leq n-1+\frac{n-\ell-1}{2}.
$$
If $k$ and $\ell$ has the same parity, then
$$
n-1\leq t(n,k,\ell)\leq n-1+\frac{n-\ell}{2}.
$$
\end{thm}

\section{Nordhaus-Guddum results}

To begin with, we have the following result.

\begin{pro}\label{pro3}
Let $G$ be a connected graph of order $n$. Then the following are
equivalent.

$(1)$ $rvx_3(G)=0$;

$(2)$ $sdiam_3(G)=2$;

$(3)$ $n-2\leq \delta(G)\leq n-1$.
\end{pro}
\begin{pf}
For Proposition \ref{pro1}, $rvx_3(G)=0$ if and only if
$sdiam_3(G)=2$. So we only need to show the equivalence of $(1)$ and
$(3)$. Suppose $n-2\leq \delta(G)\leq n-1$. Clearly, $G$ is a graph
obtained from the complete graph of order $n$ by deleting some
independent edges. For any $S=\{u,v,w\}\subseteq V(G)$, at least two
elements in $\{uv,vw,uw\}$ belong to $E(G)$. Without loss of
generality, let $uv,vw\in E(G)$. Then the tree $T$ induced by the
edges in $\{uv,vw\}$ is an $S$-Steiner tree and hence $d_G(S)\leq
2$. From the arbitrariness of $S$, we have $sdiam_3(G)\leq 2$ and
hence $sdiam_3(G)=2$. Therefore, $rvx_3(G)=0$.

Conversely, we assume $rvx_3(G)=0$. If $\delta(G)\leq n-3$, then
there exists a vertex $u\in V(G)$ such that $d_{G}(u)\leq n-3$.
Furthermore, there are two vertices, say $v,w$, such that
$uv,uw\notin E(G)$. Choose $S=\{u,v,w\}$. Clearly, any rainbow
$S$-tree must occupy at least a vertex in $V(G)\setminus S$, which
implies that $rvx_3(G)\geq 1$, a contradiction. So $n-2\leq
\delta(G)\leq n-1$.\qed \vspace{3pt}
\end{pf}

After the above preparation, we can derive a lower bound of
$rvx_3(G)+rvx_3(\overline{G})$.

\begin{lem}\label{lem1}
Let $G$ and $\overline{G}$ be a nontrivial connected graph of order
$n$. For $n\geq 5$, we have $rvx_3(G)+rvx_3(\overline{G})\geq 2$.
Moreover, the bound is sharp.
\end{lem}
\begin{pf}
From Proposition \ref{pro2}, we have $rvx_3(G)\geq 0$ and
$rvx_3(\overline{G})\geq 0$. If $rvx_3(G)=0$, then we have $n-2\leq
\delta(G)\leq n-1$ by Proposition \ref{pro3} and hence
$\overline{G}$ is disconnected, a contradiction. Similarly, we can
get another contradiction for $rvx_3(\overline{G})=0$. Therefore,
$rvx_3(G)\geq 1$ and $rvx_3(\overline{G})\geq 1$. So
$rvx_3(G)+rvx_3(\overline{G})\geq 2$.\qed\vspace{3pt}
\end{pf}

To show the sharpness of the above lower bound, we consider the
following example.\vspace{3pt}

\noindent\textbf{Example 1:} Let $H$ be a graph of order $n-4$, and
let $P=a,b,c,d$ be a path. Let $G$ be the graph obtained from $H$
and the path by adding edges between the vertex $a$ and all vertices
of $H$ and adding edges between the vertex $d$ and all vertices of
$H$; see Figure 2 $(a)$. We now show that
$rvx_3(G)=rvx_3(\overline{G})=1$. Choose $S=\{a,b,d\}$. Then any
$S$-Steiner tree must occupy at least one vertex in $V(G)\setminus
S$. Note that the vertices of $V(G)\setminus S$ in the tree must
receive different colors. Therefore, $rvx_3(G)\geq 1$. We give each
vertex in $G$ with one color and need to show that $rvx_3(G)\leq 1$.
It suffices to prove that there exists a vertex-rainbow $S$-tree for
any $S\subseteq V(G)$ with $|S|=3$. Suppose $|S\cap V(H)|=3$.
Without loss of generality, let $S=\{x,y,z\}$. Then the tree $T$
induced by the edges in $\{xa, ya,za\}$ is a vertex-rainbow
$S$-tree. Suppose $|S\cap V(H)|=2$. Without loss of generality, let
$x,y\in S\cap V(H)$. If $a\in S$, then the tree $T$ induced by the
edges in $\{xa,ya\}$ is a vertex-rainbow $S$-tree. If $b\in S$, then
the tree $T$ induced by the edges in $\{xa,ya,ab\}$ is a
vertex-rainbow $S$-tree.
\begin{figure}[!hbpt]
\begin{center}
\includegraphics[scale=0.85]{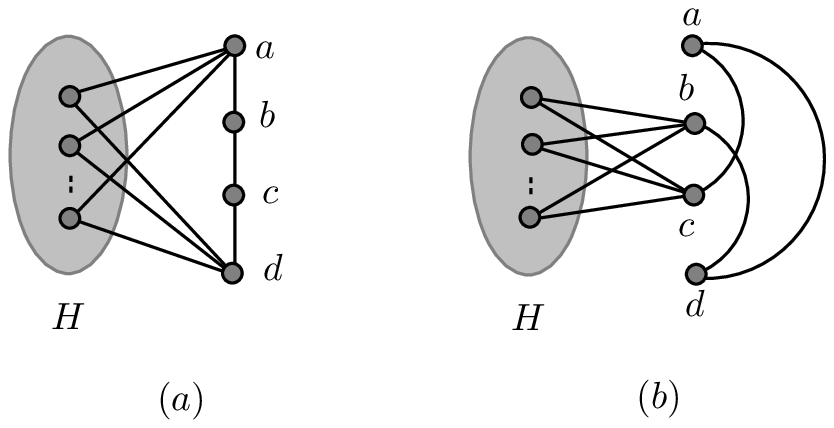}\\
Figure 2: Graphs for Example $1$.
\end{center}
\end{figure}
Suppose $|S\cap V(H)|=1$. Without loss of generality, let $x\in
S\cap V(H)$. If $a,b\in S$, then the tree $T$ induced by the edges
in $\{xa,ab\}$ is a vertex-rainbow $S$-tree. If $b,c\in S$, then the
tree $T$ induced by the edges in $\{xd,cd,bc\}$ is a vertex-rainbow
$S$-tree. If $a,c\in S$, then the tree $T$ induced by the edges in
$\{xa,ab,bc\}$ is a vertex-rainbow $S$-tree. Suppose $|S\cap
V(G')|=0$. If $a,b,c\in S$, then the tree $T$ induced by the edges
in $\{ab,bc\}$ is a vertex-rainbow $S$-tree. If $a,b,d\in S$, then
the tree $T$ induced by the edges in $\{ab,bc,cd\}$ is a
vertex-rainbow $S$-tree. From the arbitrariness of $S$, we conclude
that $rvx_3(G)\leq 1$. Similarly, one can also check that
$rvx_3(\overline{G})=1$. So $rvx_3(G)+rvx_3(\overline{G})=2$.
\qed\vspace{4pt}

We are now in a position to give an upper bound of
$rvx_3(G)+rvx_3(\overline{G})$. For $n=4$, we have
$G=\overline{G}=P_4$ since we only consider connected graphs.
Observe that $rvx_3(G)=rvx_3(\overline{G})=rvx_3(P_4)=2$.

\begin{obs}\label{obs2}
Let $G,\overline{G}$ be connected graphs of order $n \ (n=4)$. Then
$rvx_3(G)+rvx_3(\overline{G})=n$.
\end{obs}

For $n\geq 5$, we have the following upper bound of
$rvx_3(G)+rvx_3(\overline{G})$.

\begin{lem}\label{lem2}
Let $G,\overline{G}$ be connected graphs of order $n \ (n=5)$. Then
$rvx_3(G)+rvx_3(\overline{G})\leq n-1$.
\end{lem}
\begin{pf}
If $G$ is a path of order $5$, then $rvx_3(G)=3$ by Observation
\ref{obs1}. Observe that $sdiam_3(\overline{G})=3$. Then
$rvx_3(\overline{G})\leq 1$ and hence
$rvx_3(G)+rvx_3(\overline{G})\leq 4$, as desired.
\begin{figure}[!hbpt]
\begin{center}
\includegraphics[scale=0.9]{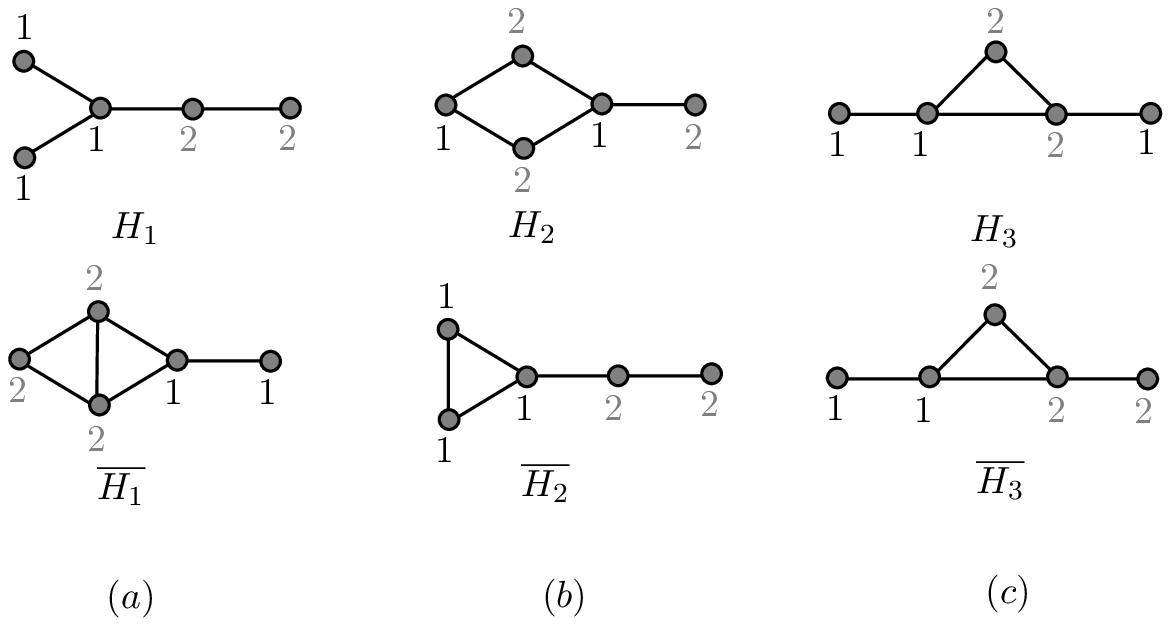}\\
Figure 3: Graphs for Lemma \ref{lem2}.
\end{center}
\end{figure}

If $G$ is a tree but not a path, then we have $G=H_1$ since
$\overline{G}$ is connected (see Figure 3 $(a)$). Clearly,
$rvx_3(G)\leq 2$. Furthermore, $\overline{G}$ consists of a $K_2$
and a $K_3$ and two edges between them (see Figure 3 $(a)$). So we
assign color $1$ to the vertices of $K_2$ and color $2$ to the
vertices of $K_3$, and this vertex-coloring makes the graph $G$
vertex-rainbow $3$-tree-connected, that is, $rvx_3(\overline{G})\leq
2$. Therefore, $rvx_3(G)+rvx_3(\overline{G})\leq 4$, as desired.

Suppose that both $G$ and $\overline{G}$ are not trees. Then
$e(G)\geq 5$ and $e(\overline{G})\geq 5$. Since
$e(G)+e(\overline{G})=e(K_5)=10$, it follows that
$e(G)=e(\overline{G})=5$. If $G$ contains a cycle of length $5$,
then $G=\overline{G}=C_5$ and hence
$rvx_3(G)=rvx_3(\overline{G})=2$. If $G$ contains a cycle of length
$4$, then $G=H_2$ (see Figure 3 $(b)$). Clearly,
$rvx_3(G)=rvx_3(\overline{G})=2$. If $G$ contains a cycle of length
$3$, then $G=\overline{G}=H_3$ (see Figure 3 $(c)$). One can check
that $rvx_3(G)=rvx_3(\overline{G})=2$. Therefore,
$rvx_3(G)+rvx_3(\overline{G})=4$, as desired.\qed
\end{pf}

\begin{lem}\label{lem3}
Let $G$ be a nontrivial connected graph of order $n$, and
$rvx_3(G)=\ell$. Let $G'$ be a graph obtained from $G$ by adding a
new vertex $v$ to $G$ and making $v$ be adjacent to $q$ vertices of
$G$. If $q \geq n-\ell$, then $rvx_3(G')\leq \ell$.
\end{lem}
\begin{pf}
Let $c: V(G)\rightarrow \{1,2,\cdots,\ell\}$ be a vertex-coloring of
$G$ such that $G$ is vertex-rainbow $3$-tree-connected. Let
$X=\{x_1, x_2,\cdots,x_q\}$ be the vertex set such that $vx_i\in
E(G')$. Set $V(G)\setminus X=\{y_1,y_2,\cdots, y_{n-q}\}$. We can
assume that there exist two vertices $y_{j_1},y_{j_2}$ such that
there is no vertex-rainbow tree connecting $\{v,y_{j_1},y_{j_2}\}$;
otherwise, the result holds obviously.
\begin{figure}[!hbpt]
\begin{center}
\includegraphics[scale=0.95]{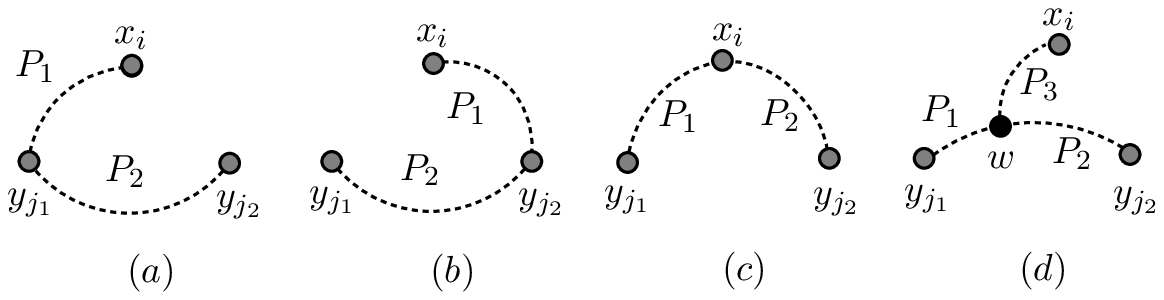}\\
Figure 4: Four type of the Steiner tree $T_i$.
\end{center}
\end{figure}

We define a minimal $S$-Steiner tree $T$ as a tree connecting $S$
whose subtree obtained by deleting any edge of $T$ does not connect
$S$. Because $G$ is vertex-rainbow $3$-tree-connected, there is a
minimal vertex-rainbow tree $T_i$ connecting
$\{x_i,y_{j_1},y_{j_2}\}$ for each $x_i \ (i\in \{1,2,\cdots,q\})$.
Then the tree $T_i$ has four types; see Figure $4$. For the type
shown in $(c)$, the Steiner tree $T_i$ connecting
$\{x_i,y_{j_1},y_{j_2}\}$ is a path induced by the edges in
$E(P_1)\cup E(P_2)$ and hence the internal vertices of the path
$T_i$ must receive different colors. Therefore, the tree induced by
the edges in $E(P_1)\cup E(P_2)\cup \{vx_i\}$ is a vertex-rainbow
tree connecting $\{v,y_{j_1},y_{j_2}\}$, a contradiction. So we only
need to consider the other three cases shown in Figure 4
$(a),(b),(d)$. Obviously, $T_i\cap T_j$ may not be empty. Then we
have the following claim.\vspace{3pt}

\noindent\textbf{Claim 1:} No other vertex in
$\{x_1,x_2,\cdots,x_q\}$ different from $x_i$ belong to $T_i$ for
each $1\leq i\leq q$. \vspace{3pt}

\noindent{\emph{Proof of Claim $1$}}: Assume, to the contrary, that
there exists a vertex $x_i'\in \{x_1,x_2,\cdots,x_q\}$ such that
$x_i'\neq x_i$ and $x_i'\in V(T_i)$. For the type shown in Figure 4
$(a)$, the Steiner tree $T_i$ connecting $\{x_i,y_{j_1},y_{j_2}\}$
is a path induced by the edges in $E(P_1)\cup E(P_2)$ and hence the
internal vertices of the path $T_i$ receive different colors. If
$x_i'\in V(P_1)$, then the tree induced by the edges in $E(P_1')\cup
E(P_2)\cup \{vx_i\}$ is a vertex-rainbow tree connecting
$\{v,y_{j_1},y_{j_2}\}$ where $P_1'$ is the path between the vertex
$x_i'$ and the vertex $y_{j_1}$ in $P_1$, a contradiction. If
$x_i'\in V(P_2)$, then the tree induced by the edges in $E(P_2)\cup
\{vx_i\}$ is a vertex-rainbow tree connecting
$\{v,y_{j_1},y_{j_2}\}$, a contradiction. The same is true for the
type shown in Figure 4 $(b)$. For the type shown in Figure 4 $(c)$,
the Steiner tree $T_i$ connecting $\{x_i,y_{j_1},y_{j_2}\}$ is a
tree induced by the edges in $E(P_1)\cup E(P_2)\cup E(P_3)$ and
hence the internal vertices of the tree $T_i$ receive different
colors. Without loss of generality, let $x_i'\in V(P_1)$. Then the
tree induced by the edges in $E(P_1')\cup E(P_2)\cup E(P_3)$ is a
vertex-rainbow tree connecting $\{v,y_{j_1},y_{j_2}\}$ where $P_1'$
is the path between the vertex $x_i'$ and the vertex $v$ in $P_1$, a
contradiction.\qed \vskip 0.5em

From Claim $1$, since there is no vertex-rainbow tree connecting
$\{v,y_{j_1},y_{j_2}\}$, it follows that there exists a vertex
$y_{k_i}$ such that $c(x_i)=c(y_{k_i})$ for each tree $T_i$, which
implies that the colors that are assigned to $X$ are among the
colors that are assigned to $V(G)\setminus X$. So $rvx_3(G)=\ell\leq
n-q$. Combining this with the hypothesis $q\geq n-\ell$, we have
$rvx_3(G)=n-q$, that is, all vertices in $V(G)\setminus X$ have
distinct colors. Now we construct a new graph $G'$, which is induced
by the edges in $E(T_1)\cup E(T_2)\cup \cdots \cup E(T_q)$.\vskip
0.5em

\noindent\textbf{Claim 2}: For every $y_t$ not in $G'$, there exists
a vertex $y_s\in G'$ such that $y_ty_s\in E(G)$.\vskip 0.5em

\noindent{\emph{Proof of Claim $2$}}: Assume, to the contrary, that
$N(y_t)\subseteq \{x_1,x_2,\cdots,x_q\}$. Since $G$ is
vertex-rainbow $3$-tree-connected, there is a vertex-rainbow tree
$T$ connecting $\{y_t,y_{j_1},y_{j_2}\}$. Let $x_r$ be the vertex in
the tree $T$ such that $x_r\in N_G(y_t)$. Then tree induced by the
edges in $(E(T)\setminus \{y_tx_r\})\cup \{vx_r\}$ is a
vertex-rainbow tree connecting $\{v,y_{j_1},y_{j_2}\}$, a
contradiction.\qed\vskip 0.5em

From Claim $2$, $G[y_1,y_2,\cdots, y_{n-q}]$ is connected. Clearly,
$G[y_1,y_2,\cdots, y_{n-q}]$ has a spanning tree $T$. Because the
tree $T$ has at least two pendant vertices, there must exist a
pendant vertex whose color is different from $x_1$, and we assign
the color to $x_1$. One can easily check that $G$ is still
vertex-rainbow $3$-tree-connected, and there is a vertex-rainbow
tree connecting $\{v,y_{j_1},y_{j_2}\}$. If there still exist two
vertices $y_{j_3},y_{j_4}$ such that there is no vertex-rainbow tree
connecting $\{v,y_{j_3},y_{j_4}\}$, then we do the same operation
until there is a vertex-rainbow tree connecting
$\{v,y_{j_r},y_{j_s}\}$ for each pair $y_{j_r},y_{j_s}\in
\{1,2,\cdots, n-q\}$. Thus $G'$ is vertex-rainbow
$3$-tree-connected. So $rvc(G')\leq \ell$. \qed
\end{pf}\vspace{5pt}

\noindent\textbf{Proof of Theorem \ref{th2}:} We prove this theorem
by induction on $n$. By Lemma \ref{lem2}, the result is evident for
$n=5$. We assume that $rvx_3(G)+rvx_3(\overline{G})\leq n-1$ holds
for complementary graphs on $n$ vertices. Observe that the union of
a connected graph $G$ and its complement $\overline{G}$ is a
complete graph of order $n$, that is, $G\cup \overline{G}=K_n$. We
add a new vertex $v$ to $G$ and add $q$ edges between $v$ and
$V(G)$. Denoted by $G'$ the resulting graph. Clearly,
$\overline{G'}$ is a graph of order $n+1$ obtained from
$\overline{G}$ by adding a new vertex $v$ to $\overline{G}$ and
adding $n-q$ edges between $v$ and $V(\overline{G})$.\vskip 0.5em

\noindent\textbf{Claim 3:} $rvx_3(G')\leq rvx_3(G)+1$ and
$rvx_3(\overline{G'})\leq rvx_3(\overline{G})+1$.\vskip 0.5em

\noindent{\itshape Proof of Claim $3$}: Let $c$ be a
$rvx_3(G)$-vertex-coloring of $G$ such that $G$ is vertex-rainbow
$3$-tree-connected. Pick up a vertex $u\in N_G(v)$ and give it a new
color. It suffices to show that for any $S\subseteq V(G')$ with
$|S|=3$, there exists a vertex-rainbow $S$-tree. If $S\subseteq
V(G)$, then there exists a vertex-rainbow $S$-tree since $G$ is
vertex-rainbow $3$-tree-connected. Suppose $S\nsubseteq V(G)$. Then
$v\in S$. Without loss of generality, let $S=\{v,x,y\}$. Since $G$
is vertex-rainbow $3$-tree-connected, there exists a vertex-rainbow
tree $T'$ connecting $\{u,x,y\}$. Then the tree $T$ induced by the
edges in $E(T')\cup \{uv\}$ is a vertex-rainbow $S$-tree. Therefore,
$rvx_3(G')\leq rvx_3(G)+1$. Similarly, $rvx_3(\overline{G'})\leq
rvx_3(\overline{G})+1$. \qed\vskip 0.5em

From Claim $3$, we have $rvx_3(G')+rvx_3(\overline{G'})\leq
rvx_3(G)+1+rvx_3(\overline{G})+1\leq n+1$. Clearly,
$rvx_3(G')+rvx_3(\overline{G'})\leq n$ except possibly when
$rvx_3(G')=rvx_3(G)+1$ and $rvx_3(\overline{G'})=
rvx_3(\overline{G})+1$. In this case, by Lemma \ref{lem3}, we have
$q\leq n-rvx_3(G)-1$ and $n-q\leq n-rvx_3(\overline{G})-1$. Thus,
$rvx_3(G)+rvx_3(\overline{G})\leq (n-1-q)+(q-1)=n-2$ and hence
$rvx_3(G')+rvx_3(\overline{G'})\leq n$, as desired. This completes
the induction. \qed \vspace{3pt}

To show the sharpness of the above bound, we consider the following
example.

\noindent\textbf{Example $2$:} Let $G$ be a path of order $n$. Then
$rvx_3(G)=n-2$. Observe that $sdiam_3(\overline{G})=3$. Then
$rvx_3(\overline{G})=1$, and so we have
$rvx_3(G)+rvx_3(\overline{G})=(n-2)+1=n-1$.

\section{The minimal size of graphs with given vertex-rainbow index}

Recall that $t(n,k,\ell)$ is the minimal size of a connected graph
$G$ of order $n$ with $rvx_k(G)\leq \ell$, where $2\leq \ell\leq
n-2$ and $2\leq k\leq n$. Let $G$ be a path of order $n$. Then
$rvx_k(G)\leq n-2$ and hence $t(n,k,n-2)\leq n-1$. Since we only
consider connected graphs, it follows that $t(n,k,n-2)\geq n-1$.
Therefore, the following result is immediate.

\begin{obs}\label{th3-1}
Let $k$ be an integer with $2\leq k\leq n$. Then
$$
t(n,k,n-2)=n-1.
$$
\end{obs}

A \emph{rose graph $R_{p}$ with $p$ petals} (or \emph{$p$-rose
graph}) is a graph obtained by taking $p$ cycles with just a vertex
in common. The common vertex is called the \emph{center} of $R_{p}$.
If the length of each cycle is exactly $q$, then this rose graph
with $p$ petals is called a \emph{$(p,q)$-rose graph}, denoted by
$R_{p,q}$. Then we have the following result.\vspace{5pt}

\noindent \textbf{Proof of Theorem \ref{th3}:} Suppose that $k$ and
$\ell$ has the different parity. Then $n-\ell-1$ is even. Let $G$ be
a graph obtained from a $(\frac{n-\ell-1}{2},3)$-rose graph
$R_{\frac{n-\ell-1}{2},3}$ and a path $P_{\ell+1}$ by identifying
the center of the rose graph and one endpoint of the path. Let $w_0$
be the center of $R_{\frac{n-\ell-1}{2},3}$, and let
$C_i=w_0v_iu_iw_0 \ (1\leq i\leq \frac{n-\ell-1}{2})$ be the cycle
of $R_{\frac{n-\ell-1}{2},3}$. Let $P_{\ell+1}=w_0w_1\cdots
w_{\ell}$ be the path of order $\ell+1$. To show the $rvx_k(G)\leq
\ell$, we define a vertex-coloring $c: V(G)\rightarrow
\{0,1,2,\cdots,\ell-1\}$ of $G$ by
$$
c(v)=\left\{
\begin{array}{ll}
i, &if~v=w_i \ (0\leq i\leq \ell-1);\\
1,&if~v=u_i~or~v=v_i \ (1\leq i\leq \frac{n-\ell-1}{2})\\
1,&if~v=w_{\ell}.
\end{array}
\right.
$$
One can easily see that there exists a vertex-rainbow $S$-tree for
any $S\subseteq V(G)$ and $|S|=3$. Therefore, $rvx_k(G)\leq \ell$
and $t(n,k,\ell)\leq n-1+\frac{n-\ell-1}{2}$.

Suppose that $k$ and $\ell$ has the same parity. Then $n-\ell$ is
even. Let $G$ be a graph obtained from a $(\frac{n-\ell}{2},3)$-rose
graph $R_{\frac{n-\ell}{2},3}$ and a path $P_{\ell}$ by identifying
the center of the rose graph and one endpoint of the path. Let $w_0$
be the center of $R_{\frac{n-\ell}{2},3}$, and let $C_i=w_0v_iu_iw_0
\ (1\leq i\leq \frac{n-\ell}{2})$ be the cycle of
$R_{\frac{n-\ell}{2},3}$. Let $P_{\ell}=w_0w_1\cdots w_{\ell-1}$ be
the path of order $\ell$. To show the $rvx_k(G)\leq \ell$, we define
a vertex-coloring $c: V(G)\rightarrow \{0,1,2,\cdots,\ell-1\}$ of
$G$ by
$$
c(v)=\left\{
\begin{array}{ll}
i, &if~v=w_i \ (0\leq i\leq \ell-1);\\
1,&if~v=u_i~or~v=v_i \ (1\leq i\leq \frac{n-\ell}{2})\\
\end{array}
\right.
$$
One can easily see that there exists a vertex-rainbow $S$-tree for
any $S\subseteq V(G)$ and $|S|=3$. Therefore, $rvx_k(G)\leq \ell$
and $t(n,k,\ell)\leq n-1+\frac{n-\ell}{2}$.\qed

\end{document}